\documentclass[aos,preprint]{imsart}

\RequirePackage[OT1]{fontenc}
\RequirePackage{amsthm,amsmath}
\RequirePackage[numbers]{natbib}
\RequirePackage[colorlinks,citecolor=blue,urlcolor=blue]{hyperref}
\usepackage{enumitem}
\usepackage{amssymb,amsmath,bm,mathrsfs,amsfonts,graphicx,multirow,stmaryrd,color}
\usepackage{times}
\usepackage{natbib,dsfont}
\usepackage{amsthm}
\usepackage{epstopdf}
\usepackage{graphics}
\usepackage{subfig}
\usepackage{float}
\usepackage{mathrsfs}
\usepackage{booktabs}
\usepackage{footnote}
\usepackage{graphicx}
\usepackage{epstopdf}
\usepackage{caption}

\arxiv{arXiv:0000.0000}

\startlocaldefs
\numberwithin{equation}{section}
\newtheorem{thm}{Theorem}[section]
\newtheorem{lem}{Lemma}[section]
\newtheorem{rmk}{Remark}[section]

\newtheorem{cor}{Corollary}[section]

\newtheorem{ass}{Assumption}

\newcommand{\diag}{\mbox{diag}}

\renewcommand{\(}{\left(}
\renewcommand{\)}{\right)}
\endlocaldefs

\begin{document}

\begin{frontmatter}
\title{Community Detection Based on the $L_\infty$ convergence of eigenvectors in DCBM}%\thanksref{T1}
\runtitle{Community Detection by SCDRE}
%\thankstext{T1}{Footnote to the title with the ``thankstext'' command.}

\begin{aug}
\author{\fnms{Yan} \snm{Liu}\thanksref{t1}\ead[label=e1]{liuy035@nenu.edu.cn}},
\author{\fnms{Zhiqiang} \snm{Hou}\thanksref{t1}\ead[label=e2]{houzq399@nenu.edu.cn}}
\author{\fnms{Zhigang} \snm{Yao}\thanksref{t3}\ead[label=e5]{zhigang.yao@nus.edu.sg}},
\author{\fnms{Zhidong} \snm{Bai}\thanksref{t1}\ead[label=e3]{baizd@nenu.edu.cn}},
\author{\fnms{Jiang} \snm{Hu}\thanksref{t2}\ead[label=e4]{huj156@nenu.edu.cn}},
\author{\fnms{Shurong} \snm{Zheng}\thanksref{t4}\ead[label=e6]{zhengsr@nenu.edu.cn}},

\thankstext{t1}{Supported by NSFC Grant 11571067}
\thankstext{t2}{Supported by NSFC Grant 11771073}
\thankstext{t3}{Supported by the MOE Tier 1 R-155-000-210-114 and Tier 2 R-155-000-184-112 at the National University of Singapore}
\thankstext{t4}{Supported by NSFC Grants 11522105 and 11690015}
\runauthor{Liu Y. et al.}

\affiliation{Northeast Normal University and National University of Singapore}

\address{KLASMOE School of Mathematics \& Statistics\\
Northeast Normal University\\
P. R. China 130024\\
\printead{e1}\\
\phantom{E-mail:\ }\printead*{e2}\\
\phantom{E-mail:\ }\printead*{e3}\\
\phantom{E-mail:\ }\printead*{e4}\\
\phantom{E-mail:\ }\printead*{e6}}

\address{Department of Statistics and Applied Probability\\
21 Lower Kent Ridge Road\\
National University of Singapore\\
Singapore 117546\\
\printead{e5}}
\end{aug}

\begin{abstract}
%Community detection is a fundamental problem in network analysis.
Spectral clustering is one of the most popular
algorithms for community detection in network analysis. Based on this rationale, in this paper we give the convergence rate of eigenvectors for the adjacency matrix in the $l_\infty$ norm,
under the stochastic block model (BM) and degree corrected stochastic block model (DCBM), adding some mild and rational conditions.  We also extend this result to a more general model, presented based on the DCBM such that the value of random variables in the adjacency matrix is not 0 or 1, but an arbitrary real number. During the process of proving the above conclusion, we obtain the relationship of the eigenvalues in the adjacency matrix and the corresponding `population' matrix, which vary in dimension from the community-wise edge probability matrix. Using that result, we can give an estimate of the number of the communities in a known set of network data. %to solve the problem that how to determine.
Meanwhile we proved the consistency of the estimator. Furthermore, according to the derivation of proof for the convergence of eigenvectors, we propose a new approach to community detection -- Spectral Clustering based on Difference of Ratios of Eigenvectors (SCDRE). Our simulation experiments demonstrate the superiority of our method
in community detection.
\end{abstract}

\begin{keyword}[class=MSC]
\kwd[Primary ]{62H30}
\kwd{91C20}
\kwd[; secondary ]{62P25}
\end{keyword}

\begin{keyword}
\kwd{community detection, spectral clustering, stochastic block model, degree corrected block model, random matrix theory, Wigner matrix.}
\end{keyword}

\end{frontmatter}

\section{Introduction}\
In recent years, the network structure has been widely used in various disciplines, such as social sciences,
computer science, physics and biology \cite{N2010,L2010,G2010,K2002}. Some network structures, for example, the biological network, scientific collaboration network, and telecommunications network, are well known and worth studying. Therefore, the corresponding network data analysis comes about due to its necessity in current times \cite{N2004,P2008,N2001}. The tools and methods for solving problems in this area play particularly important roles, especially in the presence of large datasets \cite{B2006}.

One of the fundamental problems in network analysis is community detection. The question ``What is the community'' has no precise definition, but there is an intuitive explanation connected to the graph theory. In literature such as \cite{F2010} and \cite{KL1970}, a network dataset is usually denoted as a graph (directed or undirected) $\mathscr{N} = (V ,E)$, where $V$ is a set of $n$ nodes representing the individuals to be researched and $E$ is a set of edges between the $n$ nodes, which represent the similarity (or the relationship) of the individuals. Furthermore, the ``community'' is a group of nodes that are, in some sense, more similar (more connected) to each other than the nodes in other communities. In this case, community detection in network analysis, which is somewhat analogous to clustering in multivariate analysis, aims to reorganize the network into a community structure. That is to say, the objective of community detection is to get disparate piles of nodes from $n$ nodes. The nodes with more edges between themselves are organized into an identical pile and the nodes with fewer edges are organized into the rest of the piles in the network. The ultimate goal involves estimating and labeling each node belonging to the different piles.

Various methods have been developed to solve this problem.
%Roughly, three groups of methods can be identified.
Some traditional methods of community detection are based on graph partitioning, such as the spectral bisection in computer science \cite{PASHL1990}, the Kernighan--Lin algorithm \cite{KL1970} and hierarchical
clustering \cite{CNM2004}.
%isn't ideal for the types of real-world network data(e).
Some newer methods include spectral clustering \cite{L2007} and modularity \cite{N2006,NG2004}, which involve the optimization of some reasonable global criteria over all possible network partitions. %These are ideal for our paper.
The ideas of these methods have also been integrated in our paper.
In addition, some model-based methods rely on the fitting of a probabilistic model to a network with communities, for instance the mixture models
for directed networks \cite{NG2007}, multivariate latent variable models \cite{PASHL1990}, latent feature
models \cite{H2008} and mixed membership stochastic block models for modeling overlapping
communities \cite{ABFX2008}. Two famous examples of such models are the stochastic block model \cite{KN2011,HLL1983}, or block model for short, and the degree corrected stochastic block model \cite{KN2011}, a variation of the block model. Both models are the main foci of this article, and more details are provided in Section 2.

%In this paper, due to
The outstanding performance among the other clusters and the application in community detection of spectral clustering methods based on the Laplacian \cite{NJW2002,MV2013} captures our interest.
We find that for many decades, the mass of algorithms for networks do not prove their correctness explicitly in the theoretical sense, as the steps are only given in an algorithmic sense. Theoretical correctness has been gradually established in recent years, but when extended to high dimensional, the theory is not perfect. For example, in \cite{RCY2011}, they only prove that
%\bll{in the Frobenius norm, the product of two identical normalized graph Laplacian matrices converge to the product of corresponding two population matrices}
 the eigenspace of the normalized graph Laplacian asymptotically converge to that of the population normalized graph Laplacian restricted in the Frobenius norm. %they prove that the eigenvectors of the normalized graph Laplacian asymptotically converge to the eigenvectors of a population normalized graph Laplacian only in the Frobenius norm of the matrix.
 They could explain why spectral clustering works under the block model by combining with the Davis--Kahan theorem. However, the results cannot describe the properties of the eigenvectors of the adjacency matrix straightforwardly. Hence, in this paper, we try to prove the convergence of eigenvectors in the $l_\infty$ norm by using the results of the Wigner matrix of the random matrix theory and succeed in our attempt. As a matter of fact, by re-modeling the models in this paper, we discover the essential correlation of the eigenvalues or eigenvectors between the adjacency matrix and the corresponding population matrix, and then we propose the rate of convergence for the eigenvectors when the population matrix has multiple leading eigenvalues. Using our theorem, we also prove the relationship of two eigenvectors of the adjacency matrix, which can provide a better explanation of the application of spectral clustering under the DCBM. In addition, we present a new community detection method named SCDRE based on the results of the convergence of eigenvectors in this paper.

In the methods of community detection, the general assumption is that the true total number of communities $K$ is a known constant. Therefore, the estimation of the number of communities is of great practical and theoretical importance. Many scholars have got research efforts to estimate the number of communities \cite{BP2016,ZZ2011}. We tackle this problem in Section 4, where we provide an estimator of the true number of communities in a network similar \cite{CL2014}, as a byproduct during the process of proving the convergence of the eigenvectors. The consistency is also proven.

The organization of this paper is as follows. In Section 2, we state the three kinds of models discussed in this article and illustrate some additional conditions attached to the models. The results of convergence in the $l_\infty$ norm of the eigenvalues and eigenvectors of the adjacency matrix under different models are revealed in Section 3. In Section 4, we provide a theorem to determine the number of communities and propose a new approach called the ``SCDRE'' to facilitate community detection. In the last section, we conduct some simulation experiments to explain the benefits of our approach. All proofs of the article are in the appendix .

 \emph{Notations.} In this paper, for a $n\times1$ vector $x$ and any fixed $q > 0, \|x\|_q$ denotes the $l^q$-norm. When $q=2$, omit the index and abbreviate as $\|x\|$. Either $(x)_j$ or $x(j)$ denotes the $j$th component of the vector. As $X$ is a $m\times n$ matrix, we use $\|X\|$ \and $\|X\|_F$ to denote the spectral norm: the largest singular value of $X$  and the Frobenius norm of $X$,  respectively. $X_{ij}$ or $X(i,j)$ denotes the $i$th row $j$th column element of the matrix. Throughout the paper, $C$ denotes a generic positive constant that may vary from occurrence to occurrence.

\section{Models}
%Before presenting the results in this paper,
In this section, we introduce the model lucubrated in this paper and make some reasonable and necessary assumptions about the models and communities of the network discussed.

We express a network structure a $n$-node (undirected) graph $\mathscr{N} = (V ,E)$, where
${V} = \{1, 2, \ldots, n\}$ is the set of nodes and $E$ is the set of edges. Suppose that the $n$ nodes are
split into $K$ (disjoint) communities as follows:
\begin{align*}
V = V^{(1)}\cup\ldots\cup V^{(K)},
\end{align*}
where $K$ is a constant that we need to determine. Let $n_k=|V^{(k)}|$ be the number of nodes in $k$th community. The constant $K$ is essentially an abbreviation of $K^{(n)}$ as the number of nodes $n$ in the network changes.

 Throughout the paper, we focus on the relatively ``balanced'' communities and $K$ is an invariable constant for a fixed network even though the number of individuals studied increases incrementally. More details of the assumptions are as follows:
\begin{ass}\label{a.1}
There exists a constant $C>0$ such that $\min\limits_{1\leq k\leq K}n_k\geq C\frac{n}{K}$ for all $n$.
\end{ass}
\begin{ass}\label{a.2}
The true number of communities $K=K^{(n)}$ is an invariable constant as the network size $n$ increases. The number of nodes, $n$, can continue to increase infinitely.
\end{ass}
Assumption \ref{a.1} illustrates what we mean by the balanced communities aforementioned, which basically assumes that each community has a size at least proportional to $\frac{n}{K}$.
This assumption ensures the number of nodes in each community is not too small, and thus it guarantees the accuracy of community detection.

Assumption \ref{a.2}  states the invariance of $K$ varying with $n$, the case which we are only interested in in this paper.
The study of $K$ changing with $n$ is another problem for the future.

\subsection{Stochastic Block Model}
A stochastic block model for $n$ nodes with $K$ communities is parameterized by a membership vector
$g\in\{1, 2, \ldots , K\}^{n}$ and a symmetric community-wise edge probability matrix
$P=(p_{k_1k_2})_{1\leq k_1,k_2 \leq K}\in[0, 1]^{K\times K}$.
Let $A$ be the $n\times n$ adjacency matrix of the graph $\mathscr{N}$. In the stochastic block model (BM),
 we assume that:
\begin{ass}\label{a.3}
$A$ is symmetric, with ones on the diagonal (so there are self-connections).
\end{ass}
\begin{ass}\label{a.4}
The elements of the upper triangular $\{A(i, j ) : 1 \leq i <j \leq n\}$ are independent
Bernoulli random variables satisfying
\begin{align*}
\mathbb P(A(i,j)=1)=p_{k_1k_2}, \qquad\text{if}\quad g(i)=k_1\quad\text{and}\quad g(j)=k_2,
\end{align*}
where $1\leq k_1,k_2 \leq K$, $p_{k_1k_2}$ is the $(k_1,k_2)$-element of the community-wise edge probability matrix
$P$, and $g(i)$ is the $i$th component of the membership vector.
\end{ass}
Assumption \ref{a.3} and \ref{a.4} are assumptions on the adjacency matrix $A$ in the BM.  Actually, $A(i,j)$ is equal to 0 or 1 depending on whether there is an edge between the $i$th and $j$th node. In other words, this depends on whether there is some kind of correlation (relationship) between the $i$th and $j$th individual studied. If it exists, $A(i,j)$ is equal to 1; otherwise, it is equal to 0, and the symmetry of $A$ is also obvious. Moreover, the probability of an existing connection between the $i$th and $j$th node depends on the probability of an existing connection between the communities that $i$ and $j$ belong to, which is the meaning of $p_{k_1k_2}$. The values of $g(i)$ and $k_1$ denote the community label of $i$th nodes. There are similar assumptions in other two models to be described below, and we will only explain the differences.

We shall study the convergence of eigenvectors for the adjacency matrix $A$ and the corresponding population matrix. We
regard $\mathbb{E}(A)$, the expectation of the $A$, as the population matrix.
Let $|\lambda_1|\geq\ldots\geq|\lambda_K|$ be the $K$ non-zero eigenvalues (can be negative) of $\mathbb{E}(A)$, and assume all the $n\times1$ vectors $q_1,\ldots,q_K,q_{K+1},\ldots,q_n$ are $n$ normalized eigenvectors. For $1\leq i\leq n$, let $l_i$ and the $n\times1$ vector $\eta_i$ be the eigenvalues and the corresponding normalized eigenvectors of the adjacency matrix $A$, respectively. Thus, the spectral decomposition of $\mathbb{E}(A)$ is
\begin{equation}\label{e2.1}
\mathbb{E}(A)=Q\Lambda Q^{'}=
%(q_1,\ldots,q_n){\left(
%\begin{array}{cccccc}
%\lambda_1 &&&&&\\
%     &\ddots&&&&\\
%&&\lambda_K&&&\\
% &&&0&&\\
% &&&&\ddots&\\
%&&&&&0\\
%\end{array}\right)}(q_1,\ldots,q_n)^{'}=
\sum\limits_{i=1}^K\lambda_iq_iq_i^{'},
\end{equation}
where $\Lambda_{n\times n}=\diag(\lambda_1,\ldots,\lambda_K,\underbrace {0,\ldots 0}_{n-K})$ is a diagonal matrix made up of the eigenvalues of $\mathbb{E}(A)$, and $Q=(q_1,\ldots,q_n)$ is the corresponding eigenvectors matrix.
% where $q_i$ is a $n\times1$ normalized vector corresponding to the $i$ largest eigenvalue $\lambda_i$.

%In order to facilitate study the convergence between the eigenvectors of $A$ and the corresponding 'population' matrix,
% we represent the expectation of the adjacency matrix $A$ under BM is
Besides the above spectral decomposition of $\mathbb{E}(A)$, we also find that $\mathbb{E}(A)$ can be expressed in the form of a $K\times K$ community-wise edge probability matrix $P$ as follows:
\begin{align}\label{e2.2}
\mathbb{E}(A)=\(\mathbb P\(A\(i,j\)=1\)\)_{n\times n}=FPF^{'}=\sum\limits_{k_1,k_2=1}^Kp_{k_1k_2}f_{k_1}f_{k_2}^{'},
\end{align}
where the $n\times K$ matrix $F=(f_1,\dots,f_K)$ and $f_k$ is a $n\times 1$ indicator vector of the $k$th community. The component of $f_k$ is
\begin{equation*}
f_k(i)=\left\{
\begin{aligned}{}
&1, \qquad &\text{if}\quad g(i)=k, \\
&0 , \qquad  &\text{otherwise}.\\
\end{aligned}
\right.
\end{equation*}
It is noteworthy that $f_k$ is orthogonal to each other for $1\leq k\leq K$.

To compare the two decomposition forms of $\mathbb{E}(A)$ in Equation (\ref{e2.1}) and Equation (\ref{e2.2}), we normalize vector $f_k$ and let
\begin{equation}\label{e.1}
\widetilde f_k=\frac{1}{\sqrt n_k}f_k, \qquad\alpha_k=\frac{n_k}{n},
\end{equation}
where $n_k$ is the number of nodes in $k$th community.
Hence, Equation (\ref{e2.2}) becomes
\begin{equation}\label{e2.3}
\mathbb{E}(A)=n\sum\limits_{k_1,k_2=1}^K\widetilde{p}_{k_1k_2}\widetilde f_{k_1}\widetilde f^{'}_{k_2},
\qquad\text{where}\quad\widetilde{p}_{k_1k_2}=p_{k_1k_2}\sqrt{\alpha_{k_1}\alpha_{k_2}}.
\end{equation}
Let ${\widetilde{P}}=(\widetilde{p}_{k_1k_2})_{1 \leq k_1, k_2\leq K}$. We break down ${\widetilde{P}}$ by
the spectral decomposition
\begin{equation}
{\widetilde{P}} =\sum\limits_{i=1}^K\widetilde\lambda_i\widetilde q_i\widetilde q_i^{'},
\end{equation}
and plug it into Equation (\ref{e2.3}), where $\widetilde\lambda_i$ and the $K\times1$ vector $\widetilde q_i,1\leq i\leq K$ are the eigenvalues and the
corresponding normalized eigenvectors of ${\widetilde{P}}$ respectively.
Thus, Equation (\ref{e2.2}) becomes
\begin{equation}\label{e2.5}
\mathbb{E}(A)=\sum\limits_{i=1}^Kn\widetilde\lambda_i\left(\sum\limits_{k_1=1}^K\widetilde q_i(k_1)\widetilde f_{k_1}\right)
{\left(\sum\limits_{k_2=1}^K\widetilde q_i(k_2)\widetilde f_{k_2}\right)}^{'},
\end{equation}
where $\widetilde q_i(k_1)$ is $k_1$th component of $\widetilde q_i$.

Contrasting Equations (\ref{e2.5}) and (\ref{e2.1}), it is easy to verify that for $1 \leq i\leq K$, the eigenvalues and  corresponding eigenvectors between $\mathbb{E}(A)$ and ${\widetilde{P}}$ undergo transformations as follows:
\begin{equation}\label{e2.6}
\lambda_i=n\widetilde\lambda_i,\qquad
q_i=\sum\limits_{k=1}^K\widetilde q_i(k)\widetilde f_{k}.
\end{equation}
These are important equations that run through the paper and provide clues for proving the main conclusions.

\subsection{Degree Corrected Stochastic Block Model}

A degree corrected stochastic block model (DCBM) is an extension of a block model that allows for the heterogeneity of degrees within
 a block. To be more precise, in addition to the class membership vector $g_{n\times1}$ and
community-wise edge probability matrix $P_{K\times K}$ in a stochastic block model, the DCBM has an extra degree parameter $\{\theta(i), 1\leq i\leq n\}$ called the \emph{degree heterogeneity parameter} or \emph{heterogeneity parameter}. It is necessary to add the parameter $\theta$ because the BM is inapplicable to vast real-world networks.

In the DCBM, we append the following two assumptions:
\begin{ass}\label{a.5}
 Assumption \ref{a.3} remains. In addition, the elements of the upper triangular $\{A(i,j) : 1\leq i < j\leq n\}$ are independent variables satisfying
\begin{align}\label{e..}
\mathbb P(A(i,j)=1)=\theta(i)\theta(j)p_{k_1k_2}, \qquad\text{if} \quad g(i)=k_1, g(j)=k_2.
\end{align}
\end{ass}
%\begin{rmk}
Assumption \ref{a.5} for the DCBM plays the same role as Assumption \ref{a.3} for the BM, as both are assumptions about the adjacency matrix. The only difference is the addition of parameter $\theta$ to the DCBM, which depends on the factors of each node individually considered in the model. Compared to the BM, the DCBM allows for degree heterogeneity and is much more realistic: for each node, it uses a free parameter to model the degree. An assumption about it is made in Assumption \ref{a.6}.
%\end{rmk}

Let $\theta$ and $\theta_k, k=1,\ldots,K$ be two $n\times 1$ vectors as follows:
\begin{equation*}
\theta=\left(\theta(1),\ldots,\theta(n)\right), \qquad\quad
\theta_k(i)=\left\{
\begin{aligned}{}
&\theta(i), \quad &\text{if}\quad g(i)=k, \\
&0 , \qquad  &\text{otherwise},\\
\end{aligned}
\right.
\end{equation*}
where $\theta_k(i)$ stands for the $i$th component of $\theta_k$. Under the DCBM, analogous to the decomposition method for $\mathbb{E}(A)$ under the BM, we need to define
\begin{equation*}\label{e.2}
\bar f_k=\frac{\theta_k}{\|\theta_k\|_2}, \qquad \beta_k=\frac{\|\theta_k\|_2^2}{\|\theta\|_2^2}.
\end{equation*}
\begin{ass}\label{a.6}
 For $1 \leq i\leq n$, $0<\theta(i)<1$. When $n\rightarrow\infty$, let $\beta_k\rightarrow C_k$ and $\frac{\|\theta\|^2}{n}\rightarrow C$, where $C$ and $C_k$ are two constants.
\end{ass}
\begin{rmk}
The constraint on $\theta(i)$ is because the left-hand side of Equation (\ref{e..}) is a probability. The constant $C_k$ may vary with $k$, and in fact, for $1\leq k\leq K$, we have $0<C_k<1$, and $\sum\limits_{k=1}^KC_k=1$.
\end{rmk}
Therefore, under the DCBM, the expectation of the adjacency matrix $A$ can be written as
\begin{align}\label{e2.22}
\begin{split}
\mathbb{E}(A)&=\sum\limits_{k_1,k_2=1}^Kp_{k_1k_2}\theta_{k_1}\theta_{k_2}^{'}\\
       &=\|\theta\|_2^2\sum\limits_{k_1,k_2=1}^K\overline{p}_{k_1k_2}\bar f_{k_1}\bar f_{k_2}^{'},
\qquad where\quad\overline{p}_{k_1k_2}=p_{k_1k_2}\sqrt{\beta_{k_1}\beta_{k_2}}.
\end{split}
\end{align}
Let $\overline{P}=(\overline{p}_{k_1k_2})_{1 \leq k_1, k_2\leq K}$, and let $\bar\lambda_1,\ldots,\bar\lambda_K$ and $\bar q_1,\ldots,\bar q_K$ be the eigenvalues and the corresponding normalized eigenvectors of ${\overline{P}}$. Hence, we can get the transformations of eigenvalues and eigenvectors between $\mathbb{E}(A)$ and ${\overline{P}}$  resembling the BM, such that
\begin{equation}\label{2.8}
\lambda_i=\|\theta\|^2\bar\lambda_i,\qquad
q_i=\sum\limits_{k=1}^K\bar q_i(k)\bar f_{k}.
\end{equation}
\begin{rmk}
As a matter of fact, the BM is the special case of the DCBM that regards the degree heterogeneity parameters of all nodes as equivalent.
\end{rmk}

\subsection{A more general model derived from the DCBM}

In the two models described in subsection, the values of the random variables in the adjacency matrix under the BM and DCBM are either 0, 1 or the product of the two corresponding heterogeneity parameters -- the Bernoulli distribution or a variant of Bernoulli distribution. Here, we develop a model that is more general than the previous two models. The distributions of the random variables in this model are not associated with Bernoulli, but relative to the community-wise edge probability matrix $P$.
\begin{ass}\label{a.7}
 The adjacency matrix $A$ also satisfies Assumption \ref{a.3}, and for $1\leq i, j\leq n$,
\begin{align*}
A(i,j)=x_{ij},
\end{align*}
where $x_{ij}$ is a random variable and is independent of each other such that the following two conditions:
\begin{itemize}
\item  $\mathbb{E}x_{ij}=\theta(i)\theta(j)p_{k_1,k_2}$, if $g(i)=k_1, g(j)=k_2$.
Here, $\theta(i)$ is a random variable that satisfies Assumption \ref{a.6} except for $0<\theta(i)<1$, and its value can be greater than one.\\
\item  There exists a constant $M$ and a random variable $X$ with a finite 4th moment such that for any $C>0$, we have $\mathbb P(|x_{ij}|>C)\leq M\mathbb P(|X|>C)$.
\end{itemize}
\end{ass}
In this more general model, we regard the element of adjacency matrix, $A(i,j)$, as a relatively free random variable. This can be an arbitrary real number with a distribution related to the node itself, and the community-wise edge probability that it corresponds to.
For Condition (2), the limit on the tail probability of $x_{ij}$ is weaker than the requirement for the existence of the 4th moment.

Under this model, we can get the same transformations of eigenvalues and eigenvectors in the form of Equation (\ref{2.8}) by using identical procedures and methods to the DCBM.
\section{Estimation of the difference of eigenvalues and $L_\infty$ convergence of eigenvectors}
Before considering the convergence of eigenvectors between the $\mathbb{E}(A)$ and $A$, it is necessary to figure out the property of the difference matrix between the two matrices. We will need the following decomposition:
\begin{align}\label{e3.1}
A=B+\mathbb{E}(A),\quad\text{where}\quad B=A-\mathbb{E}(A).
\end{align}
Notice that matrix $B$ is approximate to a generalized Wigner matrix under the three models that we talked about in Section 2: diagonal elements are constant, equal to 0 in fact; off-diagonals are independent of each other, expectations are 0 and variances are bounded.

We use the following two lemmas to study the properties of $B$. Moreover, Lemma \ref{l.1}, which is the order of the largest eigenvalue of $B$, is derived from the results of the limit of the largest eigenvalue of the Wigner matrix in random matrix theory literature from \cite{BY1988}.
\begin{lem}\label{l.1}
Let $B_n=\(b_{ij},1\leq i,j \leq n\)$ be a symmetric random matrix. Suppose the entries on the diagonal are constant, and the entries off the diagonal are independent such that there exists a constant $M$ and a random variable $X$ with finite 4th moment, for any constant $C>0$, we have
\begin{align*}
\mathbb P\(\left|b_{ij}\right|>C\)\leq M\mathbb P\(|X|>C\), \quad\text{for}\quad i\neq j.
\end{align*}
Then, when $n\rightarrow\infty$, we have
\begin{align}
\overline \lim\left\Vert\frac{1}{\sqrt{n}}B_n\right\Vert\leq2C, \qquad \text{a.s.},
\end{align}
where $C$ is a constant and $C=\max\limits_{i\neq j}\sqrt{Var(b_{ij})}$.
\end{lem}
\begin{rmk}
The lemma is a variant of Theorem $A$ in \cite{BY1988}. The only difference is that while the entries of matrix are not only independent but also identically distributed in \cite{BY1988}, we merely advocate independence in this lemma.
Therefore, by finding the bounds of the standard deviation of all the random variables of the matrix and modifying the proof of Theorem $A$ in \cite{BY1988} slightly, we obtain this lemma.
\end{rmk}
\begin{lem}\label{l.2}
Let $B_n=\(b_{ij},1\leq i,j \leq n\)$ be a random matrix satisfying the assumptions in Lemma \ref{l.1}.
Then, when $n\rightarrow\infty$, there exists a sequence $\varepsilon_n\rightarrow0$ such that
\begin{align*}
\mathbb P\(\cup_{i,j=1}^n|b_{ij}|\geq\varepsilon_n\sqrt{n}, i.o.\)\rightarrow0.
\end{align*}
\end{lem}
\begin{rmk}
The proof of Lemma \ref{l.2} is in the Appendix. Actually, the condition for the diagonal variables to be constant and requirement of the symmetric matrix in this lemma is not necessary, but is expressed in this way to accommodate our model assumptions.
\end{rmk}
By applying the above lemmas, we can obtain the following theorems about the convergence of eigenvectors.
\begin{thm}[$L_\infty$ convergence of eigenvectors]\label{th1}
Let $A$ be a $n\times n$ adjacency matrix of the network graph with a community structure that satisfies Assumptions \ref{a.1}, \ref{a.2}. Assume that $A$ is generated from a BM that satisfies Assumptions \ref{a.3}, \ref{a.4}, a DCBM that satisfies Assumptions \ref{a.5}, \ref{a.6}, or a more general model derived from a DCBM that satisfies Assumption \ref{a.7}.

 Let $\lambda_{k_{t-1}+1}=\ldots=\lambda_{k_t}=\gamma_t, t=1,\dots,s$ where $k_0=0, k_s=K$, $U_t=\(\eta_{k_{t-1}+1}, \dots, \eta_{k_t}\), Q_t=\(q_{k_{t-1}+1}, \dots, q_{k_t}\)$ and $T_t=Q_t^{'}U_t$. Suppose a gap between any two eigenvalues of ${P^*}$, where the generic sign $*$ stands for $\widetilde\quad$ under the BM, and $\bar\quad$ under the DCBM and the more general model, respectively. That is to say, there is the constant $C$, such that
\begin{align}\label{e3.3}
|\lambda_i^*-\lambda_j^*|\geq C, \qquad\text{for}\quad 1\leq i,j\leq K.
\end{align}
Then, the eigenvectors of $A$ almost surely converge with the eigenvectors of $\mathbb{E}(A)$ at the rate of ${n}^{-1}\log{n}$ under the $L_\infty$ norm, namely
\begin{align*}
\|U_t-Q_tT_t\|_\infty=o_{a.s.}\left({n}^{-1}\log{n}\right), \qquad\text{for}\quad t=1,\dots, s.
\end{align*}
\end{thm}
\begin{cor}\label{c1}
Under the conditions of Theorem \ref{th1}, when the eigenvalues of $\mathbb{E}(A)$  are not equal to each other, we have
\begin{align*}
\max\limits_{1\leq i\leq K}||\eta_i-q_i||_\infty=o_{a.s.}\left({n}^{-1}\log{n}\right).
\end{align*}
\end{cor}
The proof of Theorem \ref{th1} is in the Appendix. The idea of proving this theorem involves using Lemma \ref{l.1} to estimate the difference of the eigenvalues between $\mathbb{E}(A)$ and $A$, which is proven in Theorem \ref{th4.1}, and then iterating the estimation of the order of Frobenius norm to get the convergence rate of the eigenvectors satisfactorily. When the multiple of the eigenvalue is equal to 1, we can conclude that $T_i$ is a constant such that $1-T_i^2>0$ and $T_i=1+O_{a.s.}\(n^{-1}\)$ from the proof in the Theorem \ref{th1}. Then, we can easily get Corollary \ref{c1}.
\section{ Spectral Clustering based on Difference of Ratios of Eigenvectors}
In this section, we give a new method to determine the value of $K$, the true number of communities. We propose an estimator of $K$ in Theorem \ref{th4.1}, besides proving its consistency.  Furthermore, we produce an algorithm called SCDRE to facilitate community detection, and present a theorem that can provide the theoretical basis for the algorithm.

\subsection{The determination of K}

The idea that we can estimate the value of $K$ comes from the proof in Theorem \ref{th1}.
\begin{thm}\label{th4.1}
Let $A$ be a $n\times n$ adjacency matrix satisfying the conditions in Theorem \ref{th1} or Corollary \ref{c1}. Then, we have an estimator of the number of communities $K$ as follows:
\begin{align}\label{e4.1}
\hat K=\max\limits_{i}\left\{i:\left|\frac{l_i}{\bar l}\right|>C_n, \text{$C_n$ is a constant associated with $n$} \right\},
\end{align}
where $l_i$ is the $i$-th largest eigenvalue of the adjacency matrix $A$ and $\bar l$ is the average of the absolute values of all $l_i$.
Moreover, the estimator is consistent, that is:
\begin{align*}
\mathbb P(\hat K=K)\rightarrow1.
\end{align*}
\end{thm}
\begin{rmk}
The eigenvalues of $A$, $\{l_i\}_{i=1}^n$, are sorted decreasingly  by absolute value. When the absolute values are equal, the positive eigenvalue comes first.
\end{rmk}
The proof for Theorem \ref{th4.1} is in the Appendix. The choice for constant $C_n$ must satisfy these two conditions: the first is $\frac{C_n}{n}\rightarrow0$, to ensure that the leading eigenvalues are bigger than the residual eigenvalues, and the second is $C_n/\sqrt n\rightarrow\infty$, to ensure the consistency of $K$. In our simulation experiments, we choose $C_n=\delta n^{\frac{3}{4}}$, where the constant $0<\delta<1$ is a scale changing with the heterogeneity parameters, and is equal to 0.03 because of the relatively good estimator of $K$ in the simulation experiments.
%\begin{rmk}
%We can also estimate $K$ by the estimator $\hat K_0=\mathop{max}\limits_{1\leq i\leq 10}\{i:\left|\frac{l_i}{l_{i+1}}\right|\}$. The advantage is that we don't have to select $C_n$. We get around the fact that the two contiguous eigenvalues both are too small which may lead to an inaccurate estimator of $K$ by limiting the range of values of $i$. We assume that the number of communities is not very large, and no more than 10 in here. In fact, we find that the very small eigenvalues locate around $\frac{n}{2}$th in all $n$ eigenvalues, so we just have to restrict the value of $i$ to a little bit less than $\frac{n}{2}$. The reason why we don't use the above estimator is that the estimator in Th4.1 is more accurate in the simulation experiments.
%\end{rmk}
\subsection{SCDRE: A new approach to community detection}
\subsubsection{The initial idea and theoretical foundation of the algorithm}
We get the following reality through the decomposition of $\mathbb{E}(A)$, the population matrix of the adjacency matrix in Section 2. For $1\leq \zeta\leq K, 1\leq i\leq n$, the $i$th component of the $\zeta$th eigenvector of $\mathbb{E}(A)$, $q_\zeta(i)$, can be written  as follows:
\begin{align}\label{e4.2}
q_\zeta(i)=\sum\limits_{k=1}^Kq_\zeta^*(k)f_{k}^*(i),
\end{align}
where the sign $*$ has been defined in Theorem \ref{th1}, $q_\zeta^*(k)$ is the $k$th element of $q_\zeta^*$, the $\zeta$th eigenvector of ${P^*}$, and $f_{k}^*(i)$ is the $i$th element of $f_{k}^*$.

Notice that the value of $i$ can be regarded as the corresponding rank of the node. For $1\leq\zeta,\omega\leq K$ and $1\leq i\leq n$, define
\begin{align}\label{e4.3}
R_{\zeta}(i)=\frac{q_{\zeta}(i)}{\sqrt{\sum\limits_{\omega=1}^{K}q_{\omega}^2(i)}}
\end{align}
as the $i$th component of the $\zeta$th eigenvector divided by the norm of the $i$th row vector of the $n\times K$ eigenspace of the matrix $\mathbb{E}(A)$. Combine Equation (\ref{e4.2}) with (\ref{e4.3}) and notice the definition of $f_{k}^*$. Thus,
\begin{align}\label{e4.4}
R_{\zeta}(i)
%=\frac{\sum\limits_{k=1}^Kq_{\zeta}^*(k)f_{k}^*(i)}{\sqrt{\sum\limits_{\omega=1}^{K}\(\sum\limits_{k=1}^Kq_{\omega}^*(k)f_{k}^*(i)\)^2}}
=\frac{q_{\zeta}^*(k_1)}{\sqrt{\sum\limits_{\omega=1}^{K}\(q_{\omega}^*(k_1)\)^2}}
=q_{\zeta}^*(k_1),\quad\text{if}\quad g(i)=k_1,
\end{align}
where the numerator is the $k_1$th component of the $\zeta$th column vector and the denominator is the norm of the $k_1$th row vector of the eigenspace of ${P^*}$, respectively. The advantage of the definition $R_{\zeta}(i)$ is apparent from Equation (\ref{e4.4}). For any node, the
value of $R_{\zeta}(i)$ is relevant to the eigenspace of ${P^*}$. Further, it is independent of the node itself (independent of $i$ and $\theta(i)$), and only depends on the label of the community where the node is located, and an eigenvector of ${P^*}$ that we arbitrarily selected.

In addition, for $1\leq\zeta\leq K$ and $1\leq i\neq j\leq n$, $i, j$ is the coordinate of the eigenvector's component, which can also be regarded as any two of all the $n$ nodes. In fact, define
\begin{align}\label{e4.5}
RR_{\zeta}(i,j)=\frac{R_{\zeta}(i)}{R_{\zeta}(j)}
\end{align}
as the ratio of the $i$th component and the $j$th component of the $\zeta$th eigenvector of $\mathbb{E}(A)$. We get rid of some influence of the selected eigenvectors to obtain the ratio once more. Actually, we have
\begin{equation*}
RR_{\zeta}(i,j)=\left\{
\begin{aligned}{}
&1, \qquad &\text{if}\quad g(i)=g(j), \\
&C, \qquad  &\text{if}\quad g(i)\neq g(j),\\
\end{aligned}
\right.
\end{equation*}
where $C$ is a constant. Notice that $C$ may be equal to 1.

The value of $RR_{\zeta}(i,j)$ merely depends on whether both the $i$th and $j$th node are located in the same the community. In other words, if $RR_{\zeta}(i,j)\neq1$, we can declare that the two nodes are not in the same community. Thus, we get the following theorem:
\begin{thm}\label{th4.2}
Let $A$ be a $n\times n$ adjacency matrix satisfying the conditions in Theorem \ref{th1} or Corollary \ref{c1}. The corresponding population matrix $\mathbb{E}(A)$ is defined in (\ref{e2.2}) and (\ref{e2.22}). Let $RR_{\zeta}(i,j)$ be defined in (\ref{e4.5}). Then, for any two of $n$ nodes, $i$ and $j$, $1\leq i\neq j\leq n$. If they are not in an identical community, there must exist $1\leq\zeta\leq K$ such that
\begin{align*}
RR_{\zeta}(i,j)\neq 1.
\end{align*}
\end{thm}
The proof of Theorem \ref{th4.2} in the Appendix. This theorem ensures that for any two nodes that are not in the same community, there must exist one leading eigenvector of $\mathbb{E}(A)$, such that our standard formula for clustering $RR_{\zeta}(i,j)$ is not equal to 1. This theorem is the initial foundation for SCDRE, our algorithm for community detection.

Theorem \ref{th4.2} and Theorem \ref{th1} motivate us to think about $RR$ in terms of the sample. Correspondingly, we define
\begin{equation}\label{e4.6}
\widehat R_{\zeta}(i)=
\left\{
\begin{aligned}{}
&\frac{\eta_{\zeta}(i)}{\sqrt{\sum\limits_{\omega=1}^{K}\eta_{\omega}^2(i)}},  \qquad &\text{if}\quad \left|\eta_{\zeta}(i)\right|>CO_i, \\
&0, \qquad  &\text{if}\quad |\eta_{\zeta}(i)|\leq CO_i,\\
\end{aligned}
\right.
\end{equation}
where $\eta_{\zeta}(i)$ is the $i$th component of the $\zeta$th eigenvector of adjacency matrix $A$, and
\begin{equation*}
CO_i=\frac{1}{K\sqrt{n}\log{n}}\sum\limits_{\zeta=1}^K|\eta_\zeta(i)|
\end{equation*}
is a truncation to ensure that each $\hat R$ is not too small under the influence of the heterogeneity parameter, so that $\hat R$ is credible. Hence, we can get the following facts about a sample version of ${RR}_{\zeta}(i,j)$.
\begin{thm}\label{th4.3}
Let $A$ be a $n\times n$ adjacency matrix satisfying the conditions in Theorem \ref{th1} or Corollary \ref{c1}. Let the two $n\times K$ matrices made of the leading $K$ eigenvectors $U^K=(U_1, U_2, \dots,U_s)=(\mu_1^{'},\mu_2^{'},\dots,\mu_n^{'})^{'}$, $Q^K=(Q_1, Q_2, \dots,$ $Q_s)=(\upsilon_1^{'},\upsilon_2^{'},\dots,\upsilon_n^{'})^{'}$, where $\mu_i$ and $\upsilon_i, i=1,\dots,n$ are $1\times K$ vectors. Then, for $1\leq i\neq j\leq n$, we have
\begin{equation*}
\left\Vert\frac{\mu_i}{\|\mu_i\|}-\frac{\mu_j}{\|\mu_j\|}\right\Vert\left\{
\begin{aligned}{}
&\leq C\(n^{-\frac{1}{2}}\log n\), \qquad &\text{if}\quad g(i)=g(j), \\
&\geq \sqrt2-C\(n^{-\frac{1}{2}}\log n\), \quad  &\text{otherwise},\\
\end{aligned}
\right.
\end{equation*}
where $C$ is a constant.
\end{thm}
The proof is in the Appendix. Notice that $\mu_i, i=1,\dots, n$ in the above theorem is not an eigenvector, but can be regarded as the row vector of the matrix made of the first $K$ leading eigenvectors corresponding to the $i$th individual we studied. Actually, the $\zeta$th component of $\frac{\mu_i}{\|\mu_i\|}$ happens to be $\widehat R_{\zeta}(i)$. Thus, according to the above theorem, we can declare that for two individuals belonging to different communities, the norm of difference between the corresponding two normalized row vectors in the leading eigenspace of the adjacency matrix is bigger than a constant close to $\sqrt{2}$; for two individuals in identical communities, the norm has the tendency to be 0. In other words, $\widehat R_{\zeta}$ of two nodes located the same community is approximate. This result is wonderful for community detection.
\subsubsection{Steps of the algorithm}
Through the discussion in the previous subsection, the convergence of the sample version of ${RR}_{\zeta}(i,j)$ provides strong evidence that we could use the ratios of eigenvectors of the sample to detect the community structure of a network. Combined with the estimation of the true number of a network, we propose a new approach to detect the community -- Spectral Clustering based on Difference of Ratios of Eigenvectors (SCDRE). Specifically, we can the split the $n$ nodes into $\widehat K$ (disjoint) communities using the following steps:
\begin{itemize}
\item  Calculate all eigenvalues and eigenvectors of the adjacency matrix $A$. In order of magnitude from large to small using the absolute of the eigenvalue, let $|l_1|\geq\ldots\geq |l_n|$ and the corresponding unit eigenvector be $\eta_1,\ldots,\eta_n$.
\item  Use the estimator in (\ref{e4.1}) to estimate the true number of communities, let $\hat K$.
\item  For all $1\leq i\leq n$ and $1\leq\zeta\leq \hat K$, calculate $\hat R_\zeta(i)$ defined in (\ref{e4.6}) and combine it into a $n\times\hat K$ matrix. Let $\Pi$, which is a variant of the eigenspace of $A$, be made of $\eta_1,\ldots,\eta_{\hat K}$.
\item  Let $1\times K$ vector $r_i, i=1,\ldots,n$ be the $n$ row vectors of $\Pi$. Regard $r_i$ as a node in $\mathbb{R^{\hat K}}$, clustering the $r_i$ by the k-means method. The $i$ node belongs to the cluster of $r_i$.
\end{itemize}

In our algorithm, the vector $r_i$ is the truncated version of $\frac{\mu_i}{\|\mu_i\|}$. The optimization function is applied in the last step by the k-mean method. As a matter of fact, the optimization function of SCDRE can be written as:
\begin{align*}
\hat g=\arg\min\limits_{\hat g}\sum_{k=1}^{\hat K}\sum_{\hat g(j)=k}\|r_j-\overline r^{k}\|^2,
\end{align*}
where $\overline r^{k}$ denotes the average of all $r_j$ belonging to the $k$th community and $\hat g$ is the membership vector  estimated by SCDRE.

\section{Simulations}
 In this section, we do several simulations to show the superiority of our work in the article.
For each simulation experiment, the number of repetitions is 100 and we conducted the following initial setup:
\begin{enumerate}
\item Choose parameters including $n$ and $K$, the number of nodes and the true number of communities.
\item Set a $K\times K$ community-wise edge probability matrix ${P}$. In general, we set the diagonal elements of ${P}$ to be greater than the non-diagonal elements, which means that the probability of association of the nodes within an identical cluster is greater than the probability of the association of nodes between two different clusters.
\item For each node, set degree heterogeneity parameters $\theta(i)\in\(0,1\)$ and give it an original setting label $g(i)\in\{1,\dots,K\}$. In our experiment, we get the $n$ true labels and the membership vector $g_{1\times n}$ by generating $n$ uniform random integers with value of $1,\dots,K$.
\end{enumerate}
 \subsection{Simulations for the estimator of the true number of communities}

 We investigate the performance of Formula (\ref{e4.1}), the estimator of the true number of communities that we proposed in Section 4.

\textbf{Experiment 1.} Assume the heterogeneity parameter of each node generated by uniform distribution from 0.15 to 1. We investigate the estimator when the true number of communities is equal to 2 and 3, and the total number of nodes is equal to 400, 500, 1000, and 2000, respectively. Let the elements of the community-wise edge probability matrix be $P(1,1)=P(2,2)=1, P(1,2)=0.5$ when $K=2$, and add an additional setting $P(1,3)=P(2,3)=0.04$ if $K=3$. We choose $0.03n^{\frac{3}{4}}$ to be the value of $C_n$ in our estimator. The results are tabulated in Table \ref{t.*}.
\begin{table}[!htbp]
\centering
\caption{Experiment 1. In each cell, the two numbers denote the mean of the estimator of $K$ and corresponding standard deviation in the bracket.}\label{t.*}
\begin{tabular}{ccccc}
\toprule
$n$ & $n=400$& $n=500$ &$n=1000$ &$n=2000$ \\
\midrule
$K=2$& 2(0)& 2(0)& 2(0)& 2(0)\\
$K=3$& 3(0)& 3(0)& 3(0)& 3(0)\\
\bottomrule
\end{tabular}
\end{table}

We can see that over 100 repetitions, all estimator obtained from Formula (\ref{e4.1}) are correct. That is to say, the Formula (\ref{e4.1}) is a fairly good estimator of $K$, provided that the number of nodes ($n$) is large enough and $C_n$ in (\ref{e4.1}) is appropriate chosen.
 \subsection{Simulations of algorithm}
Here we compare our algorithm to other four well-known algorithms: oPCA, the most primitive spectral clustering on the adjacency matrix \cite{CFT2012}; nPCA, spectral clustering on the Laplacian matrix of the graph \cite{C1997}; SCORE, proposed in \cite{J2015}; and SCORE+ in \cite{JTL2018}. Please refer to the corresponding literature for the specific steps of the algorithms.

The simulation experiments are carried out based on the relative error rate. This means that we focus on the ratio of numbers between the wrongly estimated labels and total nodes. The number of errors considered is the minimum number of wrongly estimated labels in the sense of permutation. The details of experiment are shown below.

 \textbf{Experiment 2.} In this experiment, we investigate how SCDRE, SCORE, SCORE+, oPCA,
and nPCA perform with the classical stochastic Block Model (BM).
We set $n=1000$ and $K=2$ or $3$. Let the community-wise edge probability be the same as that of Experiment 1. The heterogeneity parameters of all nodes are equal to 0.2. The performance of the five algorithms is reflected in Table \ref{t.1}.
\begin{table}[!htbp]
\centering
\caption{Experiment 2. In each cell, the two numbers denote the mean of the relative error rate and corresponding standard deviation in the bracket.}\label{t.1}
\begin{tabular}{cccccc}
\toprule
methods & SCDRE& SCORE &SCORE+ &oPCA &nPCA \\
\midrule
$K=2$& 0.058(0.009)& 0.058(0.009)&0.056(0.009)&0.059(0.009)&0.056(0.009)\\
$K=3$& 0.114(0.049)& 0.109(0.050)&0.048(0.008)&0.121(0.052)&0.133(0.0546)\\
\bottomrule
\end{tabular}
\end{table}

\textbf{Experiment 3.} In this experiment, we investigate the performance of the five algorithms in the degree corrected stochastic block model (DCBM). We set the $n, K$ and $K\times K$ community-wise edge probability in a similar manner as in Experiment 1, and the heterogeneity parameter of the $i$th node is $\theta(i)=0.015+0.785\times\(\frac{i}{N}\)^2$. The results are shown in Table \ref{t.2}.
\begin{table}[!htbp]
\centering
\caption{Experiment 3. In each cell, the two numbers denote the mean of the relative error rate and corresponding standard deviation in the bracket.}\label{t.2}
\begin{tabular}{cccccc}
\toprule
methods & SCDRE& SCORE &SCORE+ &oPCA &nPCA \\
\midrule
$K=2$& 0.073(0.008)& 0.074(0.008)&0.076(\bf{0.060})&0.269(0.023)&0.160(0.175)\\
$K=3$& 0.096(0.009)& 0.105(0.009)&0.139(\bf{0.157})&0.385(0.012)&0.342(0.115)\\
\bottomrule
\end{tabular}
\end{table}

\textbf{Experiment 4.} In this experiment, we investigate how the different heterogeneity parameters affect the performance of the five algorithms.
We set $n=1000, K=3$ and the $K\times K$ community-wise edge probability in a similar manner as in Experiment 1. We take three kinds of heterogeneity parameters of each node. The results correspond to the three rows from the top to bottom of Table \ref{t.3}, respectively. (1) $\theta(i)=c_0+\(d_0-c_0\)\times\(\frac{i}{n}\)$, (2) $\theta(i)=c_0+\(d_0-c_0\)\times\(\frac{i}{n}\)^2$, (3) $\theta(i)=c_0I\{i\leq\frac{n}{2}\}+d_0I\{i>\frac{n}{2}\}$, where $c_0=0.015,d_0=0.8$.
\begin{table}[!htbp]
\centering
\caption{Experiment 4. In each cell, the two numbers denote the mean of the relative error rate and corresponding standard deviation in the bracket.}\label{t.3}
\begin{tabular}{cccccc}
\toprule
methods & SCDRE& SCORE &SCORE+ &oPCA &nPCA \\
\midrule
mean(SD)& 0.021(0.005)& 0.021(0.005)&0.024(\bf{0.063})&0.140(0.013)&0.021(0.005)\\
        & 0.097(0.009)& 0.106(0.009)&0.122(\bf{0.130})&0.380(0.012)&0.344(0.101)\\
        & 0.153(0.011)& 0.159(0.111)&\bf{0.208}(\bf{0.163})&0.345(0.006)&0.516(0.067)\\
\bottomrule
\end{tabular}
\end{table}

The results in Table \ref{t.1} and \ref{t.2} show that our method is not very different from SCORE and SCORE+, but slightly outperforms both of the models. Combining with Table \ref{t.3}, we observe that SCORE+ is sensitive to the heterogeneity parameters. The large variation of heterogeneity parameters of different individuals in the same network dataset could largely affect the stability of SCORE+. It is worth mentioning that a phenomenon occurs uncoincidentally in several group experiments over 100 repetitions: for the same group experiment, there have been a few repetitions that give huge errors in SCORE+. More specifically, in a few repetitions among 100 of a group experiment, the number of incorrectly estimated labels account for about a half, which overall causes the mean of repetitions large, and therefore leads to a large variance. This phenomenon happens for several different group experiments (See bold figure in Tables \ref{t.2} and \ref{t.3}). The reason for this phenomenon could be from the fact that SCORE+ is good for the network structure with weak signals, see \cite{JTL2018}. In contrast, for the current simulated data, the performance of SCDRE has always been stable. The similar performance of SCDRE will be also demonstrated in the real data experiment. Comparing with the oPCA and nPCA, the advantage of the SCDRE in the DCBM is obvious.

\textbf{Experiment 5.} In this experiment, we compare the performance of five algorithms in real data. We select eight famous real data sets collected from \url{http://www-personal.umich.edu/~mejn/netdata/} and \cite{TK2011,TM2012}. We do the same pre-processing as in \cite{JTL2018} to satisfy the premise of mutual exclusion between communities. For all the datasets, the true labels were proposed by the  authors in each corresponding literature, and here we simply use that information directly. The results are shown in Table \ref{t.5}.
\begin{table}[!htbp]
\centering
\caption{Experiment 5. In each cell, the numerator and denominator of each fraction denote the number of errors and the total number of the individuals studied in real data, respectively.}\label{t.5}
\begin{tabular}{c|ccccc}
\toprule
Dataset       & SCDRE    & SCORE    &SCORE+   &oPCA     &nPCA \\
\hline
Simmons       & 253/1137 & 268/1137 &127/1137    &442/1137    &278/1137\\
Caltech       & 178/590 & 180/590   &100/590     &221/590     &174/590 \\
Dolphins      & 1/62    &0 /62      &2/62      &12/62      &   0/62\\
Polbooks      &3/92     &1/92       &2/92      &4/92      &2/92\\
UKfaculty     &6/79     & 1/79      & 2/79     & 5/79     & 1/79\\
Football      &4/110    & 5/110     & 6/110    & 6/110    & 6/110\\
Political blogs&64/1222  & 58/1222   & 51/1222  & 437/1222 & 590/1222\\
Karate        &0/34     &0/34       &1/34      &0/34      &1/34\\
\bottomrule
\end{tabular}
\end{table}

From the results of the eight real datasets, we can see that the overall performance of SCDRE and SCORE is quite similar, whether on the big datasets (i.e., Simmons, Political blogs) or small datasets  (i.e., Karate, Dolphins). Meanwhile, SCORE+ outperforms on Simmons and Caltech, probably due to that the two datasets are weak signals. However, SCORE+ is not as good as SCDRE on Dolphins, Karate or Football. oPCA and nPCA have poor performance on Simmons and Political blogs. For all datasets, SCDRE has the minimum incorrectly estimated number on Football and Karate, and it has certain advantages on other datasets, which secures SCDRE a competitive method among the rests.

\subsection{Discussion on the hamster social network}
\begin{figure}[!htbp]
\includegraphics[width=2.5 in]{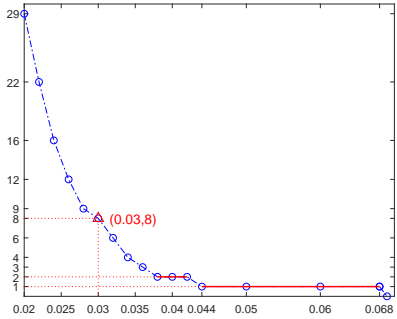}\
\captionsetup{font={footnotesize}}
\caption{\textit{This figure shows the estimation of $K$ varying with $\delta$. The the horizontal axis represents the parameter $\delta$ in $C_n$ and the vertical axis represents the the different estimation of $K$ by Formula (\ref{e4.1}).}}\label{f.1}
%\caption{table title} %%%% \scriptsize
\end{figure}

We collect the hamster social network data from \url{http://konect.cc/categories/Social/}. This data set contains friendships and family links between users of  \url{hamsterster.com}. We have no true number of the communities or the true labels of the individuals. All we know are that there are 2426 individuals in the study and the adjacency matrix established by the relationship among individuals. In the analysis, we primarily estimate the number of communities and  the label of each individual. According to Formula (\ref{e4.1}), we could obtain different estimates of $K$, by choosing different values of $\delta$ in $C_n$, where $C_n=\delta N^{\frac{3}{4}}$. The results are shown in Figure \ref{f.1}.

 In Figure \ref{f.1}, the specific range of $\delta$ has been shown since the corresponding estimation of $K$ beyond this range is almost meaningless (equals to 0 or over 30). The corresponding estimated numbers of communities $K$ for the different $\delta$ values have been expressed as hollow points. The $\delta$ is chosen from the interval [0.02,0.044] with step size of 0.002. The blue curve displays the trend of the $K$ with respect to $\delta$. We can be comparatively certain that the true number of communities should be equal to 1, 2 or 8 from Figure \ref{f.1}, which can also be obtained by observation the sequence of eigenvalues complying with Formula (\ref{e4.1}) of the adjacency matrix of the real data. We decide to ignore the case $K=1$ which means there is no community structure in the data. While other choices of $K$ are possible, we choose to only investigate the case $K=8$ since the corresponding $\delta$ happens to be 0.03 as used in Experiment 1.

 \begin{figure}[!htbp]
\includegraphics[width=2.5in]{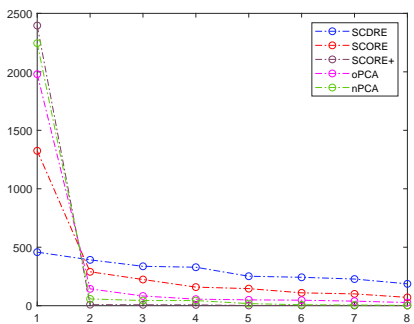}
\captionsetup{font={footnotesize}}
\caption{\textit{This figure shows the dispersion of individuals under five clustering algorithms. The horizontal axis represents the label of each community under the five algorithms, ranking from 1(maximum) to 8(minimum) by the number of individuals it contains, respectively. The vertical axis represents the number of individuals included in each community. So we have five descending curves.}}\label{f.2}
\end{figure}

 Based on this basic assumption ($K=8$), we perform the five algorithms mentioned above to estimate the label of each individual for this data. Figure \ref{f.2} shows some finding of the estimated labels. First, we obtain a relatively balanced community structure from SCDRE that the largest community contains 458 individuals, while the smallest contains 187. The number of the largest community found by SCDRE is strikingly close to the one found by the Edge Between-ness algorithms such as Girvan--Newman algorithm \cite{NG2004}. Second, the four algorithms, except SCDRE, all produce a super giant community. The largest community originate from the algorithm SCORE+, consists of 2,396 individuals, while the other seven communities have fewer than 10 individuals each (i.e., only one individual in two communities). Under oPCA and nPCA, the second largest communities are relatively prominent contain 144 and 58 individuals, respectively, although both of them are nowhere near the largest community. To some extent, the estimated labels by SCORE are the closest to SCDRE because of a relatively balanced number of individuals of the other seven communities, ranging from 297 to 71, except the largest community.

\appendix

\section{Proofs}
\begin{proof}[\text{Proof of Lemma \ref{l.2}}]
Without loss of generality, we assume $\mathbb{E}b_{ij}=0$. Otherwise, we can apply a process of centralization and it will not affect our proof. Let us fix $\varepsilon_n$ as an arbitrarily small positive constant $\varepsilon$. We have
\begin{align*}
&\mathbb P\(\cup_{i,j=1}^n|b_{ij}|\geq\varepsilon\sqrt{n}, i.o.\)
%=\mathbb P\(\bigcap_{K=1}^{\infty}\bigcup_{n=K}^{\infty}\bigcup_{i,j=1}^n|b_{ij}|\geq\varepsilon\sqrt{n}\)\\
=\lim\limits_{K\rightarrow\infty}\mathbb P\(\cup_{n=K}^{\infty}\cup_{i,j=1}^n|b_{ij}|\geq\varepsilon\sqrt{n}\)\\
\leq&\!\lim\limits_{K\rightarrow\infty}\!\sum\limits_{k=K}^{\infty}\!
\mathbb P\!\(\cup_{n=2^{k-1}}^{2^k}\!\cup_{i,j=1}^{n}\!|b_{ij}|\!\geq\varepsilon\sqrt{n}\)
\leq\!\lim\limits_{K\rightarrow\infty}\!\sum\limits_{k=K}^{\infty}\!
\mathbb P\!\(\cup_{i,j=1}^{2^k}\!|b_{ij}|\!\geq\varepsilon2^{\frac{k-1}{2}}\).
\end{align*}
Consider the series behind the limit, and according the conditions in Lemma \ref{l.2}, we have
\begin{align*}
&\sum\limits_{k=1}^{\infty}\mathbb P\(\cup_{i,j=1}^{2^k}|b_{ij}|\geq\varepsilon2^{\frac{k-1}{2}}\)\\
%&\leq C_1\sum\limits_{k=1}^{\infty}2^{2k}\mathbb P\(|X|\geq\varepsilon2^{\frac{k-1}{2}}\)\\
\leq&C_1\!\sum\limits_{k=1}^{\infty}\!2^{2k}\!\sum\limits_{m=k}^{\infty}\!
\mathbb P\!\(\varepsilon2\!^{\frac{m-1}{2}}\!\leq\!|X|\!<\!\varepsilon2\!^{\frac{m}{2}}\)\!
=\!C_1\!\sum\limits_{m=1}^{\infty}\!\mathbb P\(\varepsilon2^{\frac{m-1}{2}}\!\leq\!|X|\!<\!\varepsilon2\!^{\frac{m}{2}}\!\)
\!\sum\limits_{k=1}^{m}2^{2k}\\
%\leq& C_1\sum\limits_{m=1}^{\infty}2^{2(m+1)}\mathbb P\(\varepsilon2^{\frac{m-1}{2}}\leq|X|<\varepsilon2^{\frac{m}{2}}\)\\
\leq &4C_1\!\sum\limits_{m=1}^{\infty}2^{2m}\mathbb P\!\(\!\varepsilon2\!^{\frac{m-1}{2}}\!\leq\!|X|\!<\!\varepsilon2^{\frac{m}{2}}\)\!
\leq\!16C_1\sum\limits_{m=1}^{\infty}
\varepsilon^{-4}\mathbb{E}|X|^4I_{\{\varepsilon2^{\frac{m-1}{2}}\leq|X|<\varepsilon2^{\frac{m}{2}}\}}\\
\leq& 16C_1\sum\limits_{m=1}^{\infty}\varepsilon^{-4}\mathbb{E}|X|^4I_{\{|X|\geq\varepsilon2^{\frac{m-1}{2}}\}}<\infty.
\end{align*}
%where $X$ has finite 4th moment.
This gives $\mathbb P\(\cup_{i,j=1}^n|b_{ij}|\geq\varepsilon\sqrt{n}, i.o.\)\rightarrow0$.
Consequently, by Lemma 15 of \cite{LBH2016}, for the sequence $\varepsilon_n\rightarrow0$, we have
\begin{align*}
\mathbb P\(\cup_{i,j=1}^n|b_{ij}|\geq\varepsilon_n\sqrt{n}, i.o.\)\rightarrow0,
\end{align*}
which completes the proof.
\end{proof}
\begin{proof}[\text{Proof of Theorem \ref{th4.1}}]
We first give the reason for the expression given for the estimator of $K$, and then prove its consistency.
By the Equation (\ref{e3.1}), we have
\begin{align*}
\frac{1}{n}A=\frac{1}{n}B+\frac{1}{n}\mathbb{E}(A).
\end{align*}
For the $i$th eigenvalue of $\mathbb{E}(A)$, according to the property of eigenvalues, we have
\begin{align*}
\lambda_i=\inf\limits_{y_1,\ldots,y_{i-1}}\sup\limits_{\substack{\|x\|=1 \\ x\perp{{y_1,\ldots,y_{i-1}}}}}x^{'}\mathbb{E}(A)x.
\end{align*}
Hence, based on the triangle inequality and Lemma \ref{l.1}, we obtain
\begin{align}\label{e6.3}
\max\limits_{1\leq i\leq n}\left|\frac{1}{n}l_i-\frac{1}{n}\lambda_i\right|=\max\limits_{1\leq i\leq K}\left|\widetilde{l}_i-\widetilde\lambda_i\right|
\leq\frac{1}{n}\|B\|_2\leq\frac{C}{\sqrt{n}},
\end{align}
where $\widetilde{l}_i=\frac{1}{n}l_i$ and $\widetilde\lambda_i=\frac{1}{n}\lambda_i$ as we have defined before.
The $(n-K)$ eigenvalues of $\mathbb{E}(A)$ are 0, and lead to the first equal sign of the above formula.
 Therefore, there is a gap that order is $O\(n^\frac{1}{2}\)$ between the first $K$ dominant eigenvalues of $A$ and rest $\(n-K\)$ eigenvalues with smaller absolute value. On the account of the above reason, we can provide an estimator of $K$ in the form of Formula (\ref{e4.1}).

The consistency of the estimator remains to be proven.

 Notice that
 \begin{align*}
 \mathbb P\(\hat K= K\)=1-\mathbb P\(\hat K\neq K\)=1-\mathbb P\(\hat K>K\)-\mathbb P\(\hat K<K\).
 \end{align*}
Moreover, we have $\bar l=O\(1\)$ from the above discussion and recall the definition of $\hat K$, we have
\begin{align*}
\mathbb P\(\hat K>K\)&\leq \mathbb P\left(\underset{k>K}{\overset{n}\bigcup} \left|\frac{l_{k}}{\bar l}\right|>C_n\right)
             =\mathbb P\left(\left|\frac{l_{K+1}}{\bar l}\right|>C_n\right)
             =o(1).
\end{align*}
The first equal sign holds because of the  monotonicity of $l_i$'s . Equation (\ref{e6.3}) and the condition on $C_n$ that $\frac{C_n}{\sqrt{n}}\rightarrow\infty$ ensure $\mathbb P\left(\left|\frac{l_{K+1}}{\bar l}\right|>C_n\right)$ to be a infinitely small quantity.
On the other hand, for the same reasons as above and $\frac{C_n}{n}\rightarrow0$, we have
\begin{align*}
\mathbb P\(\hat K<K\)=
\mathbb P\(\left|\frac{l_K}{\bar l}\right|<C_n \)=o(1).
\end{align*}
The proof is complete.%Inequation \ref{e6.3} and the condition on $C_n$ that $\frac{C_n}{n}\rightarrow0$ and the definition of $\hat K$.}
\end{proof}
\begin{proof}[\text{Proof of Theorem \ref{th1}}]
We prove the result under the more general model derived from the DCBM in details below. The results in DCBM and BM come naturally using the same procedure because the Assumptions \ref{a.3}, \ref{a.4} or \ref{a.5}, \ref{a.6} are the special cases of Assumption \ref{a.7}.
%The results under the BM can be obtained by taking the equivalent value of $\theta(j), 1\leq j\leq n$ equal to 1.

Let $\widetilde{l}_i=\frac{1}{n}l_i$, for $1\leq i\leq n$ as before and $\widetilde{\gamma}_{t}=\frac{1}{n}\gamma_{t}$ , for $ 1\leq t\leq s$. By $\frac{1}{n}A\eta_i=\widetilde l_i\eta_i,1\leq i\leq n$ and Equation (\ref{e3.1}), we have
\begin{align*}
%\left(\frac{1}{n}B+\frac{1}{n}\mathbb{E}(A)\right)U_i&=\widetilde l_iU_i\\\left(\widetilde{l}_i-\widetilde\gamma_i\right)
\frac{1}{n}BU_t+\frac{1}{n}\mathbb{E}(A)U_t&=\widetilde\gamma_tU_t+U_tG_t,
\end{align*}
where the $\(k_t-k_{t-1}\)\times\(k_t-k_{t-1}\)$ matrix ${G_t}\!=\!\diag\!\left(\widetilde{l}_{k_{t-1}+1}-\widetilde\gamma_t,\!\ldots,\!\widetilde{l}_{k_t}-\widetilde\gamma_t \right)$.
Since $\frac{1}{n}\mathbb{E}(A)Q_t=\widetilde\gamma_tQ_t$, we obtain
\begin{align}\label{e6.4}
%\left(\frac{1}{n}\mathbb{E}(A)-\widetilde\gamma_iI_n\right)Ui&=G_iU_i-\frac{1}{n}BU_i\\
\left(\frac{1}{n}\mathbb{E}(A)-\widetilde\gamma_tI_n\right)\left(U_t-Q_tT_t\right)&=U_tG_t-\frac{1}{n}BU_t.
\end{align}
We first prove that
\begin{align*}
\left\Vert U_1-Q_1T_1\right\Vert^2_F=O_{a.s.}\(n^{-1}\).
\end{align*}
Replace $\frac{1}{n}\mathbb{E}(A)$ by $Q\widetilde\Lambda Q^{'}$ on the left side of Equation (\ref{e6.4}), where $\widetilde\Lambda=\diag(\widetilde\lambda_1,$ $\dots, \widetilde\lambda_K,\underbrace {0,\dots 0}_{n-K})=\diag(\widetilde\gamma_1,\ldots,\widetilde\gamma_s,\underbrace {0,\dots 0}_{n-K})$, we obtain
\begin{align*}
&\left(\frac{1}{n}\mathbb{E}(A)-\widetilde\gamma_1I_n\right)\left(U_1-Q_1T_1\right)\\
=&Q \diag\left(\bm0_{k_1},\widetilde\gamma_2-\widetilde\gamma_1, \ldots, \widetilde\gamma_s-\widetilde\gamma_1, -\widetilde\gamma_1, \ldots, -\widetilde\gamma_1\right)Q^{'}\left(U_1-Q_1T_1\right)\\
=&Q \diag\left(\bm1_{k_1},\widetilde\gamma_2-\widetilde\gamma_1, \ldots, \widetilde\gamma_s-\widetilde\gamma_1, -\widetilde\gamma_1, \ldots, -\widetilde\gamma_1\right)Q^{'}\left(U_1-Q_1T_1\right)\\
&-Q \diag\left(\bm1_{k_1},0, \ldots,0\right)Q^{'}\left(U_1-Q_1T_1\right).
\end{align*}
Notice that
\begin{align*}
Q \diag\left(\bm1_{k_1},0, \ldots,0\right)Q^{'}\left(U_1-Q_1T_1\right)=0.
\end{align*}
Therefore,
\begin{align}\label{e6.5}
U_1=Q_1T_1+\Delta,
\end{align}
where $\Delta=QDQ^{'}\left(U_1G_1-\frac{1}{n}BU_1\right)$ and $D$ is a diagonal matrix as
\begin{align*}
D=\diag\(\bm1_{k_1},\left(\widetilde\gamma_2-\widetilde\gamma_1\)^{-1}, \ldots, \left(\widetilde\gamma_s-\widetilde\gamma_1\right)^{-1}, -\widetilde\gamma_1^{-1}, \ldots, -\widetilde\gamma_1^{-1}\right).
\end{align*}
 Let $C=\max\{1,\left(\widetilde\gamma_2-\widetilde\gamma_1\right)^{-2}, \ldots, \left(\widetilde\gamma_s-\widetilde\gamma_1\right)^{-2},$ $ -\widetilde\gamma_1^{-2},\ldots, -\widetilde\gamma_1^{-2}\}$, $C$ is a finite constant because of the condition (\ref{e3.3}) and the condition on heterogeneity parameters in Assumption \ref{a.7}. % Based on Lemma 2.2 in \cite{B1995},
 Then, we have
\begin{align*}
\left\Vert U_1-Q_1T_1\right\Vert^2_F=\Vert \Delta\Vert^2_F
&\leq2C\left(\left\Vert U_1G_1\right\Vert^2_F+\left\Vert\frac{1}{n}BU_1\right\Vert^2_F\right)\\
&\leq2C\left(\frac{1}{n}\Vert U_1\Vert^2_F\right)\\
&=O_{a.s.}\left(n^{-1}\right).
\end{align*}
The last inequality is true because of the inequality in (\ref{e6.3}) and Lemma \ref{l.1}. According to
\begin{align}\label{3.1}
\left(U_1-Q_1T_1\right)^{'}\left(U_1-Q_1T_1\right)=I_{k_1}-T_1^{'}T_1=\Delta^{'}\Delta\geq0,
\end{align}
we have
\begin{align*}
\|I_{k_1}-T_1^{'}T_1\|_F=\|\Delta^{'}\Delta\|_F\leq\|\Delta\|^2_F=O_{a.s.}\(n^{-1}\).
\end{align*}

Using the same method, for $t=1,\ldots,s$, we can get
\begin{align}\label{3.2}
\left\Vert U_t-Q_tT_t\right\Vert^2_F=O_{a.s.}\left(n^{-1}\right)\quad\text{and}\quad
\|I_{k_t-k_{t-1}}-T_t^{'}T_t\|_F=O_{a.s.}\(n^{-1}\).
\end{align}

Based on Equation (\ref{3.1}), we have
\begin{align*}
%\left\Vert I_{k_1}-T_1^{'}T_1\right\Vert=\Vert G^{'}G\Vert=\Vert G\Vert^2=O_{a.s.}(n^{-1})
\Vert T_1\Vert_\infty^2&\leq k_1\Vert \diag\left(T_1^{'}T_1\right)\Vert_\infty\leq k_1\Vert \left(T_1^{'}T_1\right)\Vert_\infty\\
                       &\leq k_1\left(\Vert I_{k_1}\Vert_\infty+\Vert \Delta^{'}\Delta\Vert_F\right)\\
                       &=k_1\left(1+O_{a.s.}(n^{-1})\right).
\end{align*}

Recall that $q_i=\sum\limits_{k=1}^Kq_i^*(k)f^*_{k}$, and coupled with the condition on heterogeneity parameters in Assumption \ref{a.7}, it gives $\Vert Q_1\Vert_\infty$ $=O\left(n^{-\frac{1}{2}}\right)$. Based on Equation (\ref{e6.5}) and the triangle inequality, we have
\begin{align*}
\Vert U_1\Vert_\infty
&\leq\Vert Q_1T_1\Vert_\infty+\|\Delta\|_\infty\\
&\leq k_1\Vert Q_1\Vert_\infty\Vert T_1\Vert_\infty+\|\Delta\|_F\\
&=O\(n^{-\frac{1}{2}}\)\(k_1\left(1+O_{a.s.}(n^{-1})\)\)^{\frac{1}{2}}
+O_{a.s.}\(n^{-\frac{1}{2}}\)\\
&=O_{a.s.}\left(n^{-\frac{1}{2}}\right).
\end{align*}
Thus, we obtain
\begin{align*}
\left\Vert U_1G_1\right\Vert_\infty=O_{a.s.}\left(n^{-1}\right).
\end{align*}
We see that
\begin{align*}
\left\Vert\frac{1}{n}BU_1\right\Vert_\infty=\frac{1}{n}\max\limits_{1\leq \iota\leq n}\sum\limits_{i=1}^{k_1}\left|\sum\limits_{j=1}^{n}b_{\iota j}q_i(j)\right|,
\end{align*}
where $q_i(j)$ on behalf of the $j$th component of $q_i$ and $b_{\iota j}$ represent the $(\iota,j)$-element of $B$.
In addition, on account of the expressions for $q_i$ and $B$, it is implied that $q_i(j)=O(n^{-\frac{1}{2}})$, and
$b_{\iota j}$, the random part, $\mathbb{E}(b_{\iota j})=0,Var(b_{\iota j})=O(1)$. Furthermore, for a number $ M_n>0 $, using the Bernstein inequalities, we have
\begin{align*}
\mathbb P\left(\max\limits_{1\leq \iota\leq n}\sum\limits_{i=1}^{k_1}\left|\sum\limits_{j=1}^{n}b_{\iota j}q_i(j)\right|\geq M_n\right)
&\leq k_1\sum\limits_{\iota=1}^{n}\mathbb P\left(\left|\sum\limits_{j=1}^{n}b_{\iota j}q_i(j)\right|\geq M_n\right)\\
%&\leq k_1
&\leq\sum\limits_{\iota=1}^{n}\exp\left(\frac{-\frac{1}{2}M_n^2}{\sum\limits_{t=1}^{n}\mathbb{E}
\left(b_{\iota j}^2q_i^2(j)\right)+\frac{1}{3}\varepsilon_nM_n}\right)\\
&=o(1).
\end{align*}
We get the last line of above formula by choosing $ M_n=C\log n$, where $C>\frac{4}{3}\varepsilon_n$, and $\varepsilon_n$ is a an arbitrarily small positive constant such that  $\left|b_{\iota j}q_i(j)\right|\leq\varepsilon_n$ according to Lemma \ref{l.2}.
%According to the Lamma2, we know that, there exist an arbitrarily small positive sequence $\varepsilon_n\rightarrow0$, such that $P\(\left|b_{lt}\right|>\varepsilon_n\sqrt{n},i.o.\)\rightarrow0$.
Therefore, we will get
\begin{align*}
\left\Vert\frac{1}{n}BU_1\right\Vert_\infty=o_{a.s.}\left({n}^{-1}\log n\right).
\end{align*}
Hence, by triangle inequality, we have
\begin{align*}
\left\Vert U_1G_1-\frac{1}{n}BU_1\right\Vert_\infty=o_{a.s.}\left({n}^{-1}\log n\right).
\end{align*}
On account of $\widetilde\gamma_t=\frac{1}{n}\gamma_t=O(1)$, for $t=1,\dots, s$ and $\widetilde\gamma_1=\widetilde\lambda_1=\ldots=\widetilde\lambda_{k_1}$, we have
\begin{align*}
\Vert Q\widetilde\Lambda Q^{'}\left(U_1-Q_1T_1\right)\Vert_\infty
&=\max\limits_{1\leq h\leq n}\sum\limits_{j=1}^{k_1}\left|\sum\limits_{i=1}^{K}e_h^{'}\left(q_i\widetilde\lambda_i\right)
\left(q_i^{'}\left(U_1-Q_1T_1\right)e_j\right)\right|\\
&\leq O\left(n^{-\frac{1}{2}}\right)\|q_i\|\|U_1-Q_1T_1\|_F\\
&=O_{a.s.}\left(n^{-1}\right),
\end{align*}
where $e_i$ is a $n\times 1$ unit vector with $i$th component being 1, and others being 0.
Thus,
\begin{align*}
\|U_1-Q_1T_1\|_\infty&=\frac{1}{\left|\widetilde\gamma_i\right|}\left\Vert Q\widetilde\Lambda Q^{'}\left( U_1-Q_1T_1\right)-\left(U_1G_1-\frac{1}{n}BU_1\right)\right\Vert_\infty\\
&=o_{a.s.}\left({n}^{-1}\log n\right).
\end{align*}
Using the same method, we can get
\begin{align*}
\|U_t-Q_tT_t\|_\infty=o_{a.s.}\left({n}^{-1}\log n\right), \qquad \text{for} \quad t=1,\dots, s.
\end{align*}
The proof is finished.
\end{proof}
\begin{proof}[\text{Proof of Theorem 4.2}]
We prove this theorem by contradiction. Thus, for $1\leq i,j\leq n$, for any two nodes that are not in the same community, without loss of generality, let
\begin{align*}
g(i)=k_1\quad\text{and}\quad g(j)=k_2, \qquad \text{for}\quad  1\leq k_1\neq k_2\leq K.
\end{align*}
and for any $1\leq\zeta\leq K$, we assume that
\begin{align*}
RR_{\zeta}(i,j)=1.
\end{align*}
 We need to find a contradiction for the above assumption. So, from the above equation, we obtain
 \begin{align}\label{e6.8}
 RR_{\zeta}(i,j)=\frac{R_{\zeta}(i)}{R_{\zeta}(j)}=\frac{q_{\zeta}^*(k_1)}{q_{\zeta}^*(k_2)}=1,
 \end{align}
where $q_{\zeta}^*$ is the eigenvector of the matrix ${P^*}$.
% The orthogonality still keep for the $K$ row-vectors of the eigenvector matrix.
 Consider the $k_1$th and $k_2$th row of the eigenvector matrix of ${P^*}$. Because of the arbitrariness of $\zeta$, Equation (\ref{e6.8}) shows that all of the components for the $k_1$th row equals to that of the $k_2$th row. In other words, the $k_1$th row of the eigenvector matrix is parallel to the $k_2${th} row. It conflicts with the orthogonality of ${P^*}$. Thus, the theorem is proven.
\end{proof}
\begin{proof}[\text{Proof of Theorem \ref{th4.3}}]
Follows the notation in Theorem \ref{th1} and Theorem \ref{th4.2}. As in Theorem \ref{th1}, we only prove it under the DCBM. By Equation (\ref{3.2}), it is obvious that $T_t$ is an approximate orthogonal matrix with probability 1, and
\begin{align}\label{e6.9}
\|T_t\|=1+O_{a.s.}\(n^{-1}\),\quad t=1,\dots,s.
\end{align}
Applying Theorem \ref{th1}, we have
\begin{align}\label{p3.2.1}
 \|\mu_i-\upsilon_iT_t\|\leq\sqrt K\(n^{-1}\log n\)\quad\text{and}\quad\|\mu_i\|=\|\upsilon_i\|+o_{a.s.}\(n^{-1}\log n\),
\end{align}
for $i=1,\dots,n$ and $t$ such that $k_{t-1}+1\leq i\leq k_t$.
Notice Equation (\ref{2.8}) and the definition of $f_k^{*}$ in different models. For $1\leq i\leq n$, we have
\begin{align*}
%\upsilon_i=\upsilon_j,\quad if \quad g(i)=g(j).
\upsilon_{i}
%&=\left(\sum_{k=1}^{K} q_1^{*}(k) f_{k}^{*}(i), \sum_{k=1}^{K} q_{2}^{*}(k) f_{k}^{*}(i), \dots, \sum_{k=1}^{K} q_{K}^{*}(k) f_{k}^{*}(i)\right)\\
=f_{k}^{*}(i)\times\left(q_1^{*}(k),q_{2}^{*}(k), \ldots,q_{K}^{*}(k)\right), \quad\text{if}\quad g(i)=k,
\end{align*}
where $q_i^{*}$ is the eigenvector of ${P^*}$.
Therefore, for $1\leq i\neq j\leq n$, if $g(i)=g(j)=k_{1}$, we have
\begin{align*}
\frac{\upsilon_i}{f_{k_1}^{*}(i)}=\frac{\upsilon_j}{f_{k_1}^{*}(j)}\triangleq\widetilde\upsilon_{k_1},\qquad k_1=1,\ldots,K.
\end{align*}
$\widetilde\upsilon_{k_1}$ is a definition without loss of generality. Further, $\widetilde\upsilon_{k_1}$ is a unit vector. In other words, the vector $\frac{\upsilon_i}{f_{k_1}^{*}(i)}$ only depends on the community that the $i$th node belongs to, and is independent of $i$.
We have $\|\upsilon_i\|=O\(n^{-\frac{1}{2}}\)$ because of the expression for $q_i$ and Assumption \ref{a.6}.
Then, combining (\ref{e6.9}) and (\ref{p3.2.1}), we have
\begin{align*}
\left\Vert\frac{\mu_i}{\|\mu_i\|}-\frac{\mu_j}{\|\mu_j\|}\right\Vert
&\leq\left\Vert\frac{\mu_i}{\|\mu_i\|}-\frac{\upsilon_iT_{t_1}}{\|\upsilon_i\|}\right\Vert+
\left\Vert\frac{\mu_j}{\|\mu_j\|}-\frac{\upsilon_jT_{t_2}}{\|\upsilon_i\|}\right\Vert+
\left\Vert\frac{\upsilon_iT_{t_1}}{\|\upsilon_i\|}-\frac{\upsilon_jT_{t_2}}{\|\upsilon_i\|}\right\Vert\\
&\leq\left\Vert\frac{\widetilde\upsilon_{k_1}T_{t_1}}{\|\widetilde\upsilon_{k_1}\|}
  -\frac{\widetilde\upsilon_{k_1}T_{t_2}}{\|\widetilde\upsilon_{k_1}\|}\right\Vert
  +o_{a.s.}\(n^{-\frac{1}{2}}\log n\)\\
&=o_{a.s.}\(n^{-\frac{1}{2}}\log n\),
\end{align*}
where $t_1$ and $t_2$ are such that $k_{t_1-1}+1\leq i\leq k_{t_1}$ and $k_{t_2-1}+1\leq j\leq k_{t_2}$, respectively.

Moreover, if $g(i)=k_{1}, g(j)=k_{2}$ and $k_{1} \neq k_{2}$, we have
\begin{align*}
\|\widetilde\upsilon_i-\widetilde\upsilon_j\|^2=2,
\end{align*}
because of $\widetilde\upsilon_{k_1}\widetilde\upsilon_{k_1}^{'}=1$ and $\widetilde\upsilon_{k_1}\widetilde\upsilon_{k_2}^{'}=0$. Consequently,
we obtain
\begin{align*}
\left\Vert\frac{\mu_i}{\|\mu_i\|}-\frac{\mu_j}{\|\mu_j\|}\right\Vert
&\geq\left\Vert\frac{\widetilde\upsilon_{k_1}T_{t_1}}{\|\widetilde\upsilon_{k_1}\|}
  -\frac{\widetilde\upsilon_{k_2}T_{t_2}}{\|\widetilde\upsilon_{k_2}\|}\right\Vert
  -o_{a.s.}\(n^{-\frac{1}{2}}\log n\)\\
&\geq\sqrt2-o_{a.s.}\(n^{-\frac{1}{2}}\log n\).
\end{align*}
This completes the proof.
\end{proof}

\end{document}